\documentclass[12pt]{article}%
\usepackage{amsmath,amssymb}
\usepackage[applemac]{inputenc}
\usepackage{amsmath,amssymb,fullpage}
\usepackage{color}
\usepackage{amsmath}
\usepackage{amsfonts}
\usepackage{amssymb}
\usepackage{graphicx}%
\providecommand{\U}[1]{\protect\rule{.1in}{.1in}}

\newtheorem{definition}{Definition}[section]
\newtheorem{example}{Example}[section]
\newtheorem{theorem}[definition]{Theorem}
\newtheorem{problem}[definition]{Problem}
\newtheorem{remark}[definition]{ \it Remark}

\newtheorem{proposition}[definition]{Proposition}
\newtheorem{lemma}[definition]{Lemma}

\numberwithin{equation}{section}

\def\1B{\text{1\!\!I}}

\begin{document}

\title{ A Hida-Malliavin white noise calculus approach to optimal control}
\author{Nacira Agram$^{1,2}$ and Bernt {Ø}ksendal$^{1}$ }
\date{8 November 2018}
\maketitle

\begin{abstract}
The classical maximum principle for optimal stochastic control states that if
a control $\hat{u}$ is optimal, then the corresponding Hamiltonian has a
maximum at $u=\hat{u}$. The first proofs for this result assumed that the
control did not enter the diffusion coefficient. Moreover, it was assumed that
there were no jumps in the system. Subsequently it was discovered by Shige
Peng (still assuming no jumps) that one could also allow the diffusion
coefficient to depend on the control, provided that the corresponding adjoint
backward stochastic differential equation (BSDE) for the first order
derivative was extended to include an extra BSDE for the second order
derivatives.\newline

In this paper we present an alternative approach based on Hida-Malliavin
calculus and white noise theory. This enables us to handle the general case
with jumps, allowing both the diffusion coefficient and the jump coefficient
to depend on the control, and we do not need the extra BSDE with second order
derivatives.\newline

The result is illustrated by an example of a constrained linear-quadratic
optimal control.

\end{abstract}

\footnotetext[1]{Department of Mathematics, University of Oslo, P.O. Box 1053
Blindern, N--0316 Oslo, Norway. Email: \texttt{naciraa@math.uio.no,
oksendal@math.uio.no.}
\par
This research was carried out with support of the Norwegian Research Council,
within the research project Challenges in Stochastic Control, Information and
Applications (STOCONINF), project number 250768/F20.}

\footnotetext[2]{University of Biskra, Algeria.}



\section{Introduction}

Let $X^{u}(t)=X(t)$ be a solution of a controlled stochastic jump diffusion of
the form
\[%
\begin{cases}
dX(t)= & b(t,X(t),u(t))dt+\sigma(t,X(t),u(t))dB(t)\\
& +%
{\textstyle\int_{\mathbb{R}_{0}}}
\gamma(t,X(t),u(t),\zeta)\tilde{N}(dt,d\zeta);0\leq t\leq T,\\
X(0)= & x_{0}\in\mathbb{R}\ \text{(constant).}%
\end{cases}
\]
Here $B(t)$ and $\tilde{N}(dt,d\zeta):=N(dt,d\zeta)-\nu(d\zeta)dt$ is a
Brownian motion and an independent compensated Poisson random measure,
respectively, jointly defined on a filtered probability space $(\Omega
,\mathcal{F},\mathbb{F}=\{\mathcal{F}_{t}\}_{t\geq0},P)$ satisfying the usual
conditions. The measure $\nu$ is the Lévy measure of $N$, $T>0$ is a given
constant and $u=u(t)$ is our control process. We assume that
\[%
{\textstyle\int_{\mathbb{R}_{0}}}
\zeta^{2}\nu(d\zeta)<\infty.
\]
Now for $u$ to be admissible, we require that $u$ is $\mathbb{F}$-adapted and
that $u(t)\in V$ for all $t$ for some given Borel set $V\subset%
\mathbb{R}
$. The given coefficients $b(t,x,u)=b(t,x,u,\omega),\sigma(t,x,u)=\sigma
(t,x,u,\omega)$ and $\gamma(t,x,u,\zeta)=\gamma(t,x,u,\zeta,\omega)$ are
assumed to be $\mathbb{F}$-predictable for each given $x,u$ and $\zeta$.

\begin{problem}
We want to find $\hat{u}$ such that
\[
\underset{u\in\mathcal{A}}{\text{sup}}J(u)=J(\hat{u}),
\]
where $\mathcal{A}$ denotes the set of admissible controls, and
\[
J(u):=\mathbb{E[}%
{\textstyle\int_{0}^{T}}
f(t,X^{u}(t),u(t))dt+g(X^{u}(T))]
\]
is our performance functional, with a given $\mathbb{F}$-adapted profit rate
$f(t,x,u)=f(t,x,u,\omega)$ and a given $\mathcal{F}_{T}$-measurable terminal
payoff $g(x)=g(x,\omega)$. Such a control $\hat{u}$ (if it exists) is called
an \emph{optimal control}.\newline
\end{problem}

In the classical maximum principle for optimal control one associates to the
system a \emph{Hamiltonian} function and an \emph{adjoint BSDE}, involving the
first order derivatives of the coefficients of the system. The maximum
principle states that if $\hat{u}$ is optimal, then the corresponding
Hamiltonian has a maximum at $u=\hat{u}$. To prove this, one can perform a
so-called \emph{spike perturbation} of the optimal control, and study what
happens in the limit when the spike perturbation converges to $0$. This was
first done by Bensoussan \cite{Ben}, in the case when there are no jumps
($\gamma=0$) and when the diffusion coefficient $\sigma$ does not depend on
$u$.\newline

Subsequently it was discovered by Peng \cite{Peng} (still in the case with no
jumps) that the maximum principle could be extended to allow $\sigma$ to
depend on $u$ provided that the original adjoint BSDE was accompanied by a
second order BSDE and the Hamiltonian was extended accordingly. See e.g.
Chapter 3 in Yong and Zhou \cite{YZ} for a discussion of this.\newline

The purpose of our paper is to show that if we use spike perturbation combined
with white noise theory and the associated Hida-Malliavin calculus, we can
obtain a maximum principle similar to the classical type, with the classical
Hamiltonian and only the first order adjoint BSDE, allowing jumps and allowing
both the diffusion coefficient $\sigma$ and the jump coefficient $\gamma$ to
depend on $u$.\newline

We remark that if the set $\mathcal{A}$ of admissible control processes is
convex, we can also use convex perturbation to obtain related (albeit weaker)
versions of the maximum principle. See e.g. Bensoussan \cite{Ben} and Øksendal
and Sulem \cite{OS1} and the references therein. \newline Also note that Rong
proves in Chapter 12 in \cite{Rong} that if we have jumps in the dynamics and
the control domain is not convex, then the approach cannot allow the jump
coefficient to depend on the control.\newline

Our paper is organized as follows:

\begin{itemize}
\item In Section 2, we give a short survey of the Hida-Malliavin calculus.

\item In Section 3, we prove our main result.

\item In Section 4, we illustrate our result by an example of a constrained
linear-quadratic optimal control.
\end{itemize}

\section{A brief review of Hida-Malliavin calculus for Lévy processes}

The Malliavin derivative was originally introduced by Malliavin in \cite{M} as
a stochastic calculus of variation used to prove results about smoothness of
densities of solutions of stochastic differential equations in $\mathbb{R}%
^{n}$ driven by Brownian motion. The domain of definition of the Malliavin
derivative is a subspace $\mathbb{D}_{1,2}$ of $\mathbb{L}^{2}(P)$.
Subsequently, in Aase \textit{et al} \cite{AaOPU} the Malliavin derivative was
put into the context of the white noise theory of Hida and extended to an
operator defined on the whole of $\mathbb{L}^{2}(P)$ and with values in the
Hida space $(\mathcal{S})^{\ast}$ of stochastic distributions. This extension
is called the \emph{Hida-Malliavin} derivative. \newline

There are several advantages with working with this extended Hida-Malliavin derivative:

\begin{itemize}
\item The Hida-Malliavin derivative is defined on all of $\mathbb{L}^{2}(P)$,
and it coincides with the classical Malliavin derivative on the subspace
$\mathbb{D}_{1,2}.$

\item The Hida-Malliavin derivative combines well with the white noise
calculus, including the Skorohod integral and calculus with the Wick product
$\diamond$.

\item Moreover, it extends easily to a Hida-Malliavin derivative with respect
to a Poisson random measure.
\end{itemize}

These statements are made more precise in the following brief review, where we
recall the basic definition and properties of Hida-Malliavin calculus for Lévy
processes. The summary is partly based on Agram and Øksendal \cite{AO} and
Agram \textit{et al} \cite{AO1}, \cite{AOY2}. General references for this
presentation are Aase \textit{et al} \cite{AaOPU}, Benth \cite{B}, Lindstrøm
\textit{et al} \cite{LOU}, and the books Hida \textit{et al} \cite{HKPS} and
Di Nunno \textit{et al} \cite{DOP}.\newline

In a white noise context, the Hida-Malliavin derivative is simply a
\emph{stochastic gradient}. Equivalently, one can introduce this derivative by
means of \emph{chaos expansions}, as follows:\newline First, recall the
Lévy--Itô decomposition theorem, which states that any Lévy process $Y(t)$
with
\[
\mathbb{E}[Y^{2}(t)]<\infty,\text{ for all }t
\]
can be written
\[
Y(t)=at+bB(t)+%
{\textstyle\int_{0}^{t}}
{\textstyle\int_{\mathbb{R}_{0}}}
\zeta\tilde{N}(ds,d\zeta)
\]
with constants $a$ and $b$. In view of this we see that it suffices to deal
with Hida-Malliavin calculus for $B(\cdot)$ and for
\[
\eta(\cdot):=%
{\textstyle\int_{0}^{\cdot}}
{\textstyle\int_{\mathbb{R}_{0}}}
\zeta\tilde{N}(ds,d\zeta)
\]
separately.

\subsection{Hida-Malliavin calculus for $B(\cdot)$}

A natural starting point is the Wiener-Itô chaos expansion theorem, which
states that any $F\in\mathbb{L}^{2}(\mathcal{F}_{T},P)$ can be written
\begin{equation}
\label{eq2.1}F=%
{\textstyle\sum_{n=0}^{\infty}}
I_{n}(f_{n})
\end{equation}
for a unique sequence of symmetric deterministic functions $f_{n}\in
\mathbb{L}^{2}(\lambda^{n})$, where $\lambda$ is Lebesgue measure on $[0,T]$
and
\[
I_{n}(f_{n})=n!%
{\textstyle\int_{0}^{T}}
{\textstyle\int_{0}^{t_{n}}}
\cdots%
{\textstyle\int_{0}^{t_{2}}}
f_{n}(t_{1},\cdots,t_{n})dB(t_{1})dB(t_{2})\cdots dB(t_{n})
\]
(the $n$-times iterated integral of $f_{n}$ with respect to $B(\cdot)$) for
$n=1,2,\ldots$ and $I_{0}(f_{0})=f_{0}$ when $f_{0}$ is a constant.

Moreover, we have the isometry
\begin{equation}
\mathbb{E}[F^{2}]=||F||_{\mathbb{L}^{2}(P)}^{2}=%
{\textstyle\sum_{n=0}^{\infty}}
n!||f_{n}||_{\mathbb{L}^{2}(\lambda^{n})}^{2}. \label{eq2.2}%
\end{equation}

\begin{definition}
[Hida-Malliavin derivative $D_{t}$ with respect to $B(\cdot)$]
\hfill\break\textrm{Let $\mathbb{D}_{1,2}=\mathbb{D}_{1,2}^{(B)}$ be the space
of all $F\in\mathbb{L}^{2}(\mathcal{F}_{T},P)$ such that its chaos expansion
\eqref{eq2.1} satisfies
\[
||F||_{\mathbb{D}_{1,2}}^{2}:=%
{\textstyle\sum_{n=1}^{\infty}}
nn!||f_{n}||_{\mathbb{L}^{2}(\lambda^{n})}^{2}<\infty.
\]
}

\textrm{For $F\in\mathbb{D}_{1,2}$ and $t\in\lbrack0,T]$, we define the
\emph{Hida-Malliavin derivative} or \emph{the stochastic gradient}) of $F$ at
$t$ (with respect to $B(\cdot)$), $D_{t}F,$ by
\[
D_{t}F=%
{\textstyle\sum_{n=1}^{\infty}}
nI_{n-1}(f_{n}(\cdot,t)),
\]
where the notation $I_{n-1}(f_{n}(\cdot,t))$ means that we apply the
$(n-1)$-times iterated integral to the first $n-1$ variables $t_{1}%
,\cdots,t_{n-1}$ of $f_{n}(t_{1},t_{2},\cdots,t_{n})$ and keep the last
variable $t_{n}=t$ as a parameter.}
\end{definition}

One can easily check that
\begin{equation}
\mathbb{E}[%
{\textstyle\int_{0}^{T}}
(D_{t}F)^{2}dt]=%
{\textstyle\sum_{n=1}^{\infty}}
nn!||f_{n}||_{\mathbb{L}^{2}(\lambda^{n})}^{2}=||F||_{\mathbb{D}_{1,2}}^{2},
\label{eq2.3}%
\end{equation}
so $(t,\omega)\mapsto D_{t}F(\omega)$ belongs to $\mathbb{L}^{2}(\lambda\times
P)$.


\begin{example}
\textrm{If $F=%
{\textstyle\int_{0}^{T}}
f(t)dB(t)$ with $f\in\mathbb{L}^{2}(\lambda)$ deterministic, then
\[
D_{t}F=f(t)\mbox{ for }a.a.\,t\in\lbrack0,T].
\]
More generally, if }$\mathrm{\psi}$\textrm{$(s)$ is Itô integrable,
}$\mathrm{\psi}$\textrm{$(s)\in\mathbb{D}_{1,2}$ for $a.a.\;s$ and
$D_{t}\mathrm{\psi}(s)$ is Itô integrable for $a.a.\;t$, then
\begin{equation}
D_{t}[%
{\textstyle\int_{0}^{T}}
\mathrm{\psi}(s)dB(s)]=%
{\textstyle\int_{0}^{T}}
D_{t}\mathrm{\psi}(s)dB(s)+\mathrm{\psi}%
(t)\;\mbox{for a.a. $(t,\omega)$}.\label{eq2.4}%
\end{equation}
}
\end{example}

Some other basic properties of the Hida-Malliavin derivative $D_{t}$ are the following:

\begin{enumerate}
\item[(i)] \textbf{Chain rule } \newline Suppose $F_{1},\ldots,F_{m}%
\in\mathbb{D}_{1,2}$ and that $\Phi:\mathbb{R}^{m}\rightarrow\mathbb{R}$ is
$C^{1}$
with bounded partial derivatives. Then, $\Phi(F_{1},\cdots,F_{m})\in
\mathbb{D}_{1,2}$ and
\[
D_{t}\Phi(F_{1},\cdots,F_{m})=%
{\textstyle\sum_{i=1}^{m}}
\tfrac{\partial\Phi}{\partial x_{i}}(F_{1},\cdots,F_{m})D_{t}F_{i}.
\]

\item[(ii)] \textbf{Duality formula} \newline Suppose $\psi(t)$ is
$\mathbb{F}$-adapted with $\mathbb{E}[%
{\textstyle\int_{0}^{T}}
\psi^{2}(t)dt]<\infty$ and let $F\in\mathbb{D}_{1,2}$. Then,
\[
\mathbb{E}[F%
{\textstyle\int_{0}^{T}}
\psi(t)dB(t)]=\mathbb{E}[%
{\textstyle\int_{0}^{T}}
\psi(t)D_{t}Fdt].
\]

\item[(iii)] \textbf{Malliavin derivative and adapted processes}\newline If
$\varphi$ is an $\mathbb{F}$-adapted process, then%
\[
D_{s}\varphi(t)=0\text{ for }s>t.
\]

\end{enumerate}

\begin{remark}
We put $D_{t}\varphi(t)=\underset{s\rightarrow t-}{\lim}D_{s}\varphi(t)$ (if
the limit exists).
\end{remark}

\subsection{Extension to a white noise setting}

In the following, we let $(\mathcal{S})^{\ast}$ denote the Hida space of
stochastic distributions.\newline It was proved in Aase et al \cite{AaOPU}
that one can extend the Hida-Malliavin derivative operator $D_{t}$ from
$\mathbb{D}_{1,2}$ to all of $\mathbb{L}^{2}(\mathcal{F}_{T},P)$ in such a way
that, also denoting the extended operator by $D_{t}$, for all $F\in
\mathbb{L}^{2}(\mathcal{F}_{T},P)$, we have
\begin{equation}
D_{t}F\in(\mathcal{S})^{\ast}\text{ and }(t,\omega)\mapsto\mathbb{E}%
[D_{t}F\mid\mathcal{F}_{t}]\text{ belongs to }\mathbb{L}^{2}(\lambda\times P).
\label{eq2.10a}%
\end{equation}
Moreover, the following \emph{generalized Clark-Haussmann-Ocone formula} was
proved:
\begin{equation}
F=\mathbb{E}[F]+%
{\textstyle\int_{0}^{T}}
\mathbb{E}[D_{t}F\mid\mathcal{F}_{t}]dB(t) \label{eq2.11a}%
\end{equation}
for all $F\in\mathbb{L}^{2}(\mathcal{F}_{T},P)$. See Theorem 3.11 in Aase
\textit{et al} \cite{AaOPU} and also Theorem 6.35 in Di Nunno \textit{et al}
\cite{DOP}.\newline We can use this to get the following extension of the
duality formula (ii) above:

\begin{proposition}
[The generalized duality formula]Let $F\in\mathbb{L}^{2}(\mathcal{F}_{T},P)$
and let $\varphi(t,\omega)\in\mathbb{L}^{2}(\lambda\times P)$ be $\mathbb{F}%
$-adapted. Then
\begin{equation}
\mathbb{E}[F%
{\textstyle\int_{0}^{T}}
\varphi(t)dB(t)]=\mathbb{E}[%
{\textstyle\int_{0}^{T}}
\mathbb{E}[D_{t}F\mid\mathcal{F}_{t}]\varphi(t)dt]. \label{geduB}%
\end{equation}

\end{proposition}

\noindent{Proof.} \quad By \eqref{eq2.10a} and \eqref{eq2.11a} and the Itô
isometry, we get
\begin{align*}
&  \mathbb{E}[F%
{\textstyle\int_{0}^{T}}
\varphi(t)dB(t)]=\mathbb{E}[(\mathbb{E}[F]+%
{\textstyle\int_{0}^{T}}
\mathbb{E}[D_{t}F\mid\mathcal{F}_{t}]dB(t))(%
{\textstyle\int_{0}^{T}}
\varphi(t)dB(t))]\\
&  =\mathbb{E}[%
{\textstyle\int_{0}^{T}}
\mathbb{E}[D_{t}F\mid\mathcal{F}_{t}]\varphi(t)dt].
\end{align*}
\hfill$\square$ \bigskip\newline We will use this extension of the
Hida-Malliavin derivative from now on.\newline

\subsection{Hida-Malliavin calculus for $\tilde N(\cdot)$}

The construction of a stochastic derivative/Hida-Malliavin derivative in the
pure jump martingale case follows the same lines as in the Brownian motion
case. In this case, the corresponding Wiener-Itô Chaos Expansion Theorem
states that any $F\in\mathbb{L}^{2}(\mathcal{F}_{T},P)$ (where, in this case,
$\mathcal{F}_{t}=\mathcal{F}_{t}^{(\tilde{N})}$ is the $\sigma-$algebra
generated by $\eta(s):=%
{\textstyle\int_{0}^{s}}
{\textstyle\int_{\mathbb{R}_{0}}}
\zeta\tilde{N}(dr,d\zeta);\;0\leq s\leq t$) can be written as%

\begin{equation}
F=%
{\textstyle\sum_{n=0}^{\infty}}
I_{n}(f_{n});\;f_{n}\in\hat{L}^{2}((\lambda\times\nu)^{n}), \label{eq2.8}%
\end{equation}
where $\hat{L}^{2}((\lambda\times\nu)^{n})$ is the space of functions
$f_{n}(t_{1},\zeta_{1},\ldots,t_{n},\zeta_{n})$; $t_{i}\in\lbrack0,T]$,
$\zeta_{i}\in\mathbb{R}_{0}$ for $i=1,..,n$, such that $f_{n}\in\mathbb{L}%
^{2}((\lambda\times\nu)^{n})$ and $f_{n}$ is symmetric with respect to the
pairs of variables $(t_{1},\zeta_{1}),\ldots,(t_{n},\zeta_{n}).$

It is important to note that in this case, the $n-$times iterated integral
$I_{n}(f_{n})$ is taken with respect to $\tilde{N}(dt,d\zeta)$ and not with
respect to $d\eta(t).$ Thus, we define
\[
I_{n}(f_{n}):=n!%
{\textstyle\int_{0}^{T}}
\!\!%
{\textstyle\int_{\mathbb{R}_{0}}}
\!%
{\textstyle\int_{0}^{t_{n}}}
\!%
{\textstyle\int_{\mathbb{R}_{0}}}
\cdots%
{\textstyle\int_{0}^{t_{2}}}
\!\!%
{\textstyle\int_{\mathbb{R}_{0}}}
f_{n}(t_{1},\zeta_{1},\cdots,t_{n},\zeta_{n})\tilde{N}(dt_{1},d\zeta
_{1})\cdots\tilde{N}(dt_{n},d\zeta_{n}),
\]
for $f_{n}\in\hat{L}^{2}((\lambda\times\nu)^{n}).$

The Itô isometry for stochastic integrals with respect to $\tilde{N}%
(dt,d\zeta)$ then gives the following isometry for the chaos expansion:
\[
||F||_{\mathbb{L}^{2}(P)}^{2}=%
{\textstyle\sum_{n=0}^{\infty}}
n!||f_{n}||_{\mathbb{L}^{2}((\lambda\times\nu)^{n})}^{2}.
\]
As in the Brownian motion case, we use the chaos expansion to define the
Malliavin derivative. Note that in this case, there are two parameters
$t,\zeta,$ where $t$ represents time and $\zeta\neq0$ represents a generic
jump size.

\begin{definition}
[Hida-Malliavin derivative $D_{t,\zeta}$ with respect to $\tilde{N}%
(\cdot,\cdot)$]
Let $\mathbb{D}_{1,2}^{(\tilde{N})}$ be the space of all $F\in\mathbb{L}%
^{2}(\mathcal{F}_{T},P)$ such that its chaos expansion \eqref{eq2.8}
satisfies
\[
||F||_{\mathbb{D}_{1,2}^{(\tilde{N})}}^{2}:=%
{\textstyle\sum_{n=1}^{\infty}}
nn!||f_{n}||_{\mathbb{L}^{2}((\lambda\times\nu)^{n})}^{2}<\infty.
\]
For $F\in\mathbb{D}_{1,2}^{(\tilde{N})}$, we define the Hida-Malliavin
derivative of $F$ at $(t,\zeta)$ (with respect to $\tilde{N}(\cdot,\cdot))$,
$D_{t,\zeta}F,$ by
\[
D_{t,\zeta}F:=%
{\textstyle\sum_{n=1}^{\infty}}
nI_{n-1}(f_{n}(\cdot,t,\zeta)),
\]
where $I_{n-1}(f_{n}(\cdot,t,\zeta))$ means that we perform the $(n-1)-$times
iterated integral with respect to $\tilde{N}$ to the first $n-1$ variable
pairs $(t_{1},\zeta_{1}),\cdots,(t_{n},\zeta_{n}),$ keeping $(t_{n},\zeta
_{n})=(t,\zeta)$ as a parameter.

\end{definition}

In this case, we get the isometry.
\[
\mathbb{E}[%
{\textstyle\int_{0}^{T}}
{\textstyle\int_{\mathbb{R}_{0}}}
(D_{t,\zeta}F)^{2}\nu(d\zeta)dt]=%
{\textstyle\sum_{n=0}^{\infty}}
nn!||f_{n}||_{\mathbb{L}^{2}((\lambda\times\nu)^{n})}^{2}=||F||_{\mathbb{D}%
_{1,2}^{(\tilde{N})}}^{2}.
\]
(Compare with \eqref{eq2.3}.)

\begin{example}
\textrm{If $F=%
{\textstyle\int_{0}^{T}}
{\textstyle\int_{\mathbb{R}_{0}}}
f(t,\zeta)\tilde{N}(dt,d\zeta)$ for some deterministic $f(t,\zeta
)\in\mathbb{L}^{2}(\lambda\times\nu)$, then
\[
D_{t,\zeta}F=f(t,\zeta)\mbox{ for }a.a.\,(t,\zeta).
\]
More generally, if $\psi(s,\zeta)$ is integrable with respect to $\tilde
{N}(ds,d\zeta)$, $\psi(s,\zeta)\in\mathbb{D}_{1,2}^{(\tilde{N})}$ for
$a.a.\,s,\zeta$ and $D_{t,\zeta}\psi(s,\zeta)$ is integrable for
$a.a.\,(t,\zeta)$, then
\begin{equation}
D_{t,\zeta}(%
{\textstyle\int_{0}^{T}}
\!%
{\textstyle\int_{\mathbb{R}_{0}}}
\psi(s,\zeta)\tilde{N}(ds,d\zeta))=%
{\textstyle\int_{0}^{T}}
{\textstyle\int_{\mathbb{R}_{0}}}
D_{t,\zeta}\psi(s,\zeta)\tilde{N}(ds,d\zeta)+\psi(t,\zeta
)\;\mbox{ for }a.a.\,t,\zeta.
\end{equation}
}
\end{example}

The properties of $D_{t,\zeta}$ corresponding to those of $D_{t}$ are the following:

\begin{itemize}
\item[(i)] \textbf{Chain rule }\newline Suppose $F_{1},\cdots,F_{m}%
\in\mathbb{D}_{1,2}^{(\tilde{N})}$ and that $\phi:\mathbb{R}^{m}%
\rightarrow\mathbb{R}$ is continuous and bounded. Then, $\phi(F_{1}%
,\cdots,F_{m})\in\mathbb{D}_{1,2}^{(\tilde{N})}$ and
\begin{equation}
D_{t,\zeta}\phi(F_{1},\cdots,F_{m})=\phi(F_{1}+D_{t,\zeta}F_{1},\ldots
,F_{m}+D_{t,\zeta}F_{m})-\phi(F_{1},\ldots,F_{m}).
\end{equation}

\item[(ii)] \textbf{Duality formula}\newline Suppose $\Psi(t,\zeta)$ is
$\mathbb{F}$-adapted and $\mathbb{E}[%
{\textstyle\int_{0}^{T}}
{\textstyle\int_{\mathbb{R}_{0}}^{2}}
\Psi(t,\zeta)\nu(d\zeta)dt]<\infty$ and let $F\in\mathbb{D}_{1,2}^{(\tilde
{N})}$. Then,
\[
\mathbb{E}[F%
{\textstyle\int_{0}^{T}}
{\textstyle\int_{\mathbb{R}_{0}}}
\Psi(t,\zeta)\tilde{N}(dt,d\zeta)]=\mathbb{E}[%
{\textstyle\int_{0}^{T}}
{\textstyle\int_{\mathbb{R}_{0}}}
\Psi(t,\zeta)D_{t,\zeta}F\text{ }\nu(d\zeta)dt].
\]

\item[(iii)] \textbf{Hida-Malliavin derivative and adapted processes }\newline
If $\varphi$ is an $\mathbb{F}$-adapted process, then,
\[
D_{s,\zeta}\varphi(t)=0\text{ for all }s>t, \zeta\in\mathbb{R}_{0}.
\]

\end{itemize}

\begin{remark}
We put $D_{t,\zeta}\varphi(t)=\underset{s\rightarrow t-}{\lim}D_{s,\zeta
}\varphi(t)$ ( if the limit exists).
\end{remark}

\subsection{Extension to a white noise setting}

As in section 2.2, we note that there is an extension of the Hida-Malliavin
derivative $D_{t,\zeta}$ from $\mathbb{D}_{1,2}^{(\tilde{N})}$ to
$\mathbb{L}^{2}(\lambda\times P)$ such that the following extension of the
duality theorem holds:

\begin{proposition}
[Generalized duality formula]Suppose $\Psi(t,\zeta)$ is $\mathbb{F}$-adapted
and
\[
\mathbb{E}[%
{\textstyle\int_{0}^{T}}
{\textstyle\int_{\mathbb{R}_{0}}}
\Psi^{2}(t,\zeta)\nu(d\zeta)dt]<\infty,
\]
and let $F\in\mathbb{L}^{2}(\lambda\times P)$. Then,
\begin{equation}
\mathbb{E}[F%
{\textstyle\int_{0}^{T}}
{\textstyle\int_{\mathbb{R}_{0}}}
\Psi(t,\zeta)\tilde{N}(dt,d\zeta)]=\mathbb{E}[%
{\textstyle\int_{0}^{T}}
{\textstyle\int_{\mathbb{R}_{0}}}
\Psi(t,\zeta)\mathbb{E}[D_{t,\zeta}F\mid\mathcal{F}_{t}]\nu(d\zeta)dt].
\label{geduN}%
\end{equation}

\end{proposition}

\ \newline Accordingly, note that from now on we are working with this
generalized version of the Malliavin derivative. We emphasize that this
generalized Hida-Malliavin derivative $DX$ (where $D$ stands for $D_{t}$ or
$D_{t,\zeta}$, depending on the setting) exists for all $X\in\mathbb{L}%
^{2}(P)$ as an element of the Hida stochastic distribution space
$(\mathcal{S})^{\ast}$, and it has the property that the conditional
expectation $\mathbb{E}[DX|\mathcal{F}_{t}]$ belongs to $\mathbb{L}%
^{2}(\lambda\times P)$, where $\lambda$ is Lebesgue measure on $[0,T]$.
Therefore, when using the Hida-Malliavin derivative, combined with conditional
expectation, no assumptions on Hida-Malliavin differentiability in the
classical sense are needed; we can work on the whole space of random variables
in $\mathbb{L}^{2}(P)$.\newline

\subsection{Representation of solutions of BSDE}

The following result, due to Øksendal and Røse \cite{OR}, is crucial for our method:

\begin{theorem}
\label{Th2.9} Suppose that $f,p,q$ and $r$ are given càdlàg adapted processes
in $\mathbb{L}^{2}(\lambda\times P),\mathbb{L}^{2}(\lambda\times
P),\mathbb{L}^{2}(\lambda\times P)$ and $\mathbb{L}^{2}(\lambda\times\nu\times
P)$ respectively, and they satisfy a BSDE of the form
\begin{equation}%
\begin{cases}
dp(t) & =f(t)dt+q(t)dB(t)+%
{\textstyle\int_{\mathbb{R}_{0}}}
r(t,\zeta)\tilde{N}(dt,d\zeta);0\leq t\leq T,\\
p(T) & =F\in\mathbb{L}^{2}(\mathcal{F}_{T},P).
\end{cases}
\label{BSDE}%
\end{equation}

Then for a.a. $t$ and $\zeta$ the following holds:
\begin{equation}
q(t)=D_{t}p(t^{+}):=\underset{\varepsilon\rightarrow0^{+}}{\lim}%
D_{t}p(t+\varepsilon)\text{ (limit in }(\mathcal{S})^{\ast}), \label{eq2.24}%
\end{equation}%
\begin{equation}
q(t)=\mathbb{E}[D_{t}p(t^{+})|\mathcal{F}_{t}]:=\underset{\varepsilon
\rightarrow0^{+}}{\lim}\mathbb{E}[D_{t}p(t+\varepsilon)|\mathcal{F}_{t}]\text{
(limit in }\mathbb{L}^{2}(P)), \label{eq2.25}%
\end{equation}

and
\begin{equation}
r(t,\zeta)=D_{t^{,}\zeta}p(t^{+}):=\underset{\varepsilon\rightarrow0^{+}}%
{\lim}D_{t,\zeta}p(t+\varepsilon)\text{ (limit in }(\mathcal{S})^{\ast}),
\label{eq2.26}%
\end{equation}%
\begin{equation}
r(t,\zeta)=\mathbb{E}[D_{t^{,}\zeta}p(t^{+})|\mathcal{F}_{t}]:=\underset
{\varepsilon\rightarrow0^{+}}{\lim}\mathbb{E}[D_{t,\zeta}p(t+\varepsilon
)|\mathcal{F}_{t}]\text{ (limit in }\mathbb{L}^{2}(P)). \label{eq2.27}%
\end{equation}

\end{theorem}

\section{The spike variation stochastic maximum principle}

Throughout this work, we will use the following spaces:

\begin{itemize}
\item $\mathcal{S}^{2}$ is the set of ${\mathbb{R}}$-valued $\mathbb{F}%
$-adapted càdlàg processes $(X(t))_{t\in\lbrack0,T]}$ such that
\[
{\Vert X\Vert}_{\mathcal{S}^{2}}^{2}:={\mathbb{E}}[\sup_{t\in\lbrack
0,T]}|X(t)|^{2}]<\infty.
\]

\item $\mathbb{L}^{2}$ is the set of ${\mathbb{R}}$-valued $\mathbb{F}%
$-predictable processes $(Q(t))_{t\in\lbrack0,T]}$ such that
\[
\Vert Q\Vert_{\mathbb{L}^{2}}^{2}:={\mathbb{E}}[%
{\textstyle\int_{0}^{T}}
|Q(t)|^{2}dt]<\infty.
\]

\item $\mathbb{L}_{\nu}^{2}$ is the set of $\mathbb{F}$-predictable processes
$r:[0,T]\times\mathbb{R}_{0}\rightarrow\mathbb{R}$ such that
\[
||r||_{\mathbb{L}_{\nu}^{2}}^{2}:=\mathbb{E}[%
{\textstyle\int_{0}^{T}}
{\textstyle\int_{\mathbb{R}_{0}}}
|r(t,\zeta)|^{2}\nu(d\zeta)dt]<\infty.
\]

\item $\mathcal{A}$ is a set of all $\mathbb{F}$-predictable processes $u$
required to have values in a Borel set $V\subset%
\mathbb{R}
$. We call $\mathcal{A}$ the set of admissible control processes $u(\cdot
)$.\newline
\end{itemize}

The state of our system $X^{u}(t)=X(t)$ satisfies the following SDE%
\begin{equation}
\left\{
\begin{array}
[c]{ll}%
dX(t) & =b(t,X(t),u(t))dt+\sigma(t,X(t),u(t))dB(t)\\
& +\int_{\mathbb{R}_{0}}\gamma(t,X(t),u(t),\zeta)\tilde{N}(dt,d\zeta);0\leq
t\leq T,\\
X(0) & =x_{0}\in\mathbb{R}\ (\text{constant}),
\end{array}
\right.  \label{sde}%
\end{equation}
where $b(t,x,u)=b(t,x,u,\omega):\left[  0,T\right]  \times\mathbb{R}\times
U\times\Omega\rightarrow\mathbb{R}$, $\sigma(t,x,u)=\sigma(t,x,u,\omega
):\left[  0,T\right]  \times\mathbb{R}\times U\times\Omega\rightarrow%
\mathbb{R}
$ and $\gamma(t,x,u,\zeta)=:\left[  0,T\right]  \times\mathbb{R}\times U\times%
\mathbb{R}
_{0}\times\Omega\rightarrow%
\mathbb{R}
$.\newline From now on we fix an open convex set $U$ such that $V\subset U$
and we assume that $b$, $\sigma$ and $\gamma$ are continuously differentiable
and admits uniformly bounded partial derivatives in $U$ with respect to $x$
and $u$.\newline Moreover, we assume that the coefficients $b$, $\sigma$ and
$\gamma$ are $\mathbb{F}$-adapted, and uniformly Lipschitz continuous with
respect to $x$, in the sense that there is a constant $C$ such that, for all
$t\in\lbrack0,T],u\in V,\zeta\in\mathbb{R}_{0},\,x,x^{\prime}\in\mathbb{R}$ we
have
\[%
\begin{array}
[c]{c}%
\left\vert b\left(  t,x,u\right)  -b\left(  t,x^{\prime},u\right)  \right\vert
^{2}+\left\vert \sigma\left(  t,x,u\right)  -\sigma\left(  t,x^{\prime
},u\right)  \right\vert ^{2}\\
+%
{\textstyle\int_{\mathbb{R}_{0}}}
\left\vert \gamma\left(  t,x,u,\zeta\right)  -\gamma\left(  t,x^{\prime
},u,\zeta\right)  \right\vert ^{2}\nu(d\zeta)\leq C\left\vert x-x^{\prime
}\right\vert ^{2},\text{a.s.}%
\end{array}
\]
Under this assumption, there is a unique solution $X\in\mathcal{S}^{2}$ to the
equation $\left(  \ref{sde}\right)  $, such that%
\[%
\begin{array}
[c]{c}%
X(t)=x_{0}+{%
{\textstyle\int_{0}^{t}}
}b(s,X(s),u(s))ds+{%
{\textstyle\int_{0}^{t}}
}\sigma(s,X(s),u(s))dB(s)\\
+{%
{\textstyle\int_{0}^{t}}
}%
{\textstyle\int_{\mathbb{R}_{0}}}
\gamma(s,X(s),u(s),\zeta)\tilde{N}(ds,d\zeta);0\leq t\leq T.
\end{array}
\]
\newline The \emph{performance functional} has the form%
\begin{equation}%
\begin{array}
[c]{ll}%
J(u) & =\mathbb{E[}{%
{\textstyle\int_{0}^{T}}
}\text{ }f(t,X(t),u(t))dt+g(X(T))],\quad u\in\mathcal{A},
\end{array}
\label{perf}%
\end{equation}
with given functions $f:\left[  0,T\right]  \times\mathbb{R}\times
U\times\Omega\rightarrow\mathbb{R}$ and $g:\Omega\times\mathbb{R}%
\rightarrow\mathbb{R},$ assumed to be $\mathbb{F}$-adapted and $\mathcal{F}%
_{T}$-measurable, respectively, and continuously differentiable with respect
to $x$ and $u$ with bounded partial derivatives in $U$.\newline Suppose that
$\hat{u}$ is an optimal control. Fix $\tau\in\lbrack0,T),0<\epsilon<T-\tau$
and a bounded $\mathcal{F}_{\tau}$-measurable $v$ and define the spike
perturbed $u^{\epsilon}$ of the optimal control $\hat{u}$ by%

\begin{equation}
u^{\epsilon}(t)=%
\begin{cases}
\hat{u}(t); & t\in\lbrack0,\tau)\cup(\tau+\epsilon,T],\\
v; & t\in\lbrack\tau,\tau+\epsilon].
\end{cases}
\label{eq3.3}%
\end{equation}
Let $X^{\epsilon}(t):=X^{u^{\epsilon}}(t)$ and $\hat{X}(t):=X^{\hat{u}}(t)$ be
the solutions of $\left(  \ref{sde}\right)  $\ corresponding to $u=u^{\epsilon
}$ and $u=\hat{u}$, respectively.\newline Define%

\begin{equation}
Z^{\epsilon}(t):=X^{\epsilon}(t)-\hat{X}(t);\text{ }t\in\lbrack0,T].
\label{eq3.4}%
\end{equation}
Then by the \emph{mean value theorem} \footnote{Recall that if a function $f$
is continuously differentiable on an open convex set $U\subset\mathbb{R}^{n}$
and continuous on the closure $\bar{U}$, then for all $x,y\in\bar{U}$ there
exists a point $\tilde{x}$ on the straight line connecting $x$ and $y$ such
that
\begin{equation}
f(y)-f(x)=f^{\prime}(\tilde{x})(y-x):=\sum_{i=1}^{n}\frac{\partial f}{\partial
x_{i}}(\tilde{x})(y_{i}-x_{i})
\end{equation}
}, we can write%
\[%
\begin{array}
[c]{ll}%
b^{\epsilon}(t)-\hat{b}(t) & =\frac{\partial\tilde{b}}{\partial x}%
(t)Z^{\epsilon}(t)+\frac{\partial\tilde{b}}{\partial u}(t)(u^{\epsilon
}(t)-\hat{u}(t)),
\end{array}
\]
where
\[
b^{\epsilon}(t)=b(t,X^{\epsilon}(t),u^{\epsilon}(t)),\hat{b}(t)=b(t,\hat
{X}(t),\hat{u}(t)),
\]
and
\[%
\begin{array}
[c]{ll}%
\frac{\partial\tilde{b}}{\partial x}(t) & =\frac{\partial b}{\partial
x}(t,x,u)_{x=\tilde{X}(t),u=\tilde{u}(t)},
\end{array}
\]
and%
\[%
\begin{array}
[c]{ll}%
\frac{\partial\tilde{b}}{\partial u}(t) & =\frac{\partial b}{\partial
u}(t,x,u)_{x=\tilde{X}(t),u=\tilde{u}(t)}.
\end{array}
\]
Here $(\tilde{u}(t),\tilde{X}(t))$ is \emph{a point on the straight line}
between $(\hat{u}(t),\hat{X}(t))$ and $(u^{\epsilon}(t),X^{\epsilon}(t))$.
With a similar notation for $\sigma$ and $\gamma$, we get%

\begin{align}
Z^{\epsilon}(t)  &  ={%
{\textstyle\int_{\tau}^{t}}
}\{\tfrac{\partial\tilde{b}}{\partial x}(s)Z^{\epsilon}(s)+\tfrac
{\partial\tilde{b}}{\partial u}(s)(u^{\epsilon}(s)-\hat{u}(s))\}ds+{%
{\textstyle\int_{\tau}^{t}}
}\{\tfrac{\partial\tilde{\sigma}}{\partial x}(s)Z^{\epsilon}(s)+\tfrac
{\partial\tilde{\sigma}}{\partial u}(s)(u^{\epsilon}(s)-\hat{u}%
(s))\}dB(s)\nonumber\label{eq3.5}\\
&  +{%
{\textstyle\int_{\tau}^{t}}
}%
{\textstyle\int_{\mathbb{R}_{0}}}
\{\tfrac{\partial\tilde{\gamma}}{\partial x}(s,\zeta)Z^{\epsilon}%
(s)+\tfrac{\partial\tilde{\gamma}}{\partial u}(s,\zeta)(u^{\epsilon}%
(s)-\hat{u}(s))\}\tilde{N}(ds,d\zeta);\tau\leq t\leq\tau+\epsilon,
\end{align}
and%

\begin{equation}%
\begin{array}
[c]{l}%
Z^{\epsilon}(t)={%
{\textstyle\int_{\tau+\epsilon}^{t}}
}\tfrac{\partial\tilde{b}}{\partial x}(s)Z^{\epsilon}(s)ds+%
{\textstyle\int_{\tau+\epsilon}^{t}}
\tfrac{\partial\tilde{\sigma}}{\partial x}(s)Z^{\epsilon}(s)dB(s)\\
\text{ \ \ \ \ \ \ \ \ \ \ \ }+%
{\textstyle\int_{\tau+\epsilon}^{t}}
{\textstyle\int_{\mathbb{R}_{0}}}
\tfrac{\partial\tilde{\gamma}}{\partial x}(s,\zeta)(s)Z^{\epsilon}(s)\tilde
{N}(ds,d\zeta);\tau+\epsilon\leq t\leq T.
\end{array}
\label{eq3.6}%
\end{equation}
On other words,%
\begin{equation}
\left\{
\begin{array}
[c]{ll}%
dZ^{\epsilon}(t) & =\{\tfrac{\partial\tilde{b}}{\partial x}(t)Z^{\epsilon
}(t)+\tfrac{\partial\tilde{b}}{\partial u}(t)(v-\hat{u}(t))\}dt+\{\tfrac
{\partial\tilde{\sigma}}{\partial x}(t)Z^{\epsilon}(t)+\tfrac{\partial
\tilde{\sigma}}{\partial u}(t)(v-\hat{u}(t))\}dB(t)\\
& +%
{\textstyle\int_{\mathbb{R}_{0}}}
\{\tfrac{\partial\tilde{\gamma}}{\partial x}(t,\zeta)Z^{\epsilon}%
(t)+\tfrac{\partial\tilde{\gamma}}{\partial u}(t,\zeta)(v-\hat{u}%
(t))\}\tilde{N}(dt,d\zeta);\tau\leq t\leq\tau+\epsilon,
\end{array}
\right.  \label{var1}%
\end{equation}
and
\begin{equation}%
\begin{array}
[c]{ll}%
dZ^{\epsilon}(t) & =\tfrac{\partial\tilde{b}}{\partial x}(t)Z^{\epsilon
}(t)dt+\tfrac{\partial\tilde{\sigma}}{\partial x}(t)Z^{\epsilon}(t)dB(t)+%
{\textstyle\int_{\mathbb{R}_{0}}}
\tfrac{\partial\tilde{\gamma}}{\partial x}(t,\zeta)Z^{\epsilon}(t)\tilde
{N}(dt,d\zeta);\tau+\epsilon\leq t\leq T.
\end{array}
\label{var2}%
\end{equation}

\begin{remark}
\-

\begin{enumerate}
\item Note that since the process
\[
\eta(t):=%
{\textstyle\int_{0}^{t}}
{\textstyle\int_{\mathbb{R}_{0}}}
\zeta\tilde{N}(ds,d\zeta);t\geq0
\]
is a Lévy process, we know that for every given (deterministic) time $t\geq0$
the probability that $\eta$ jumps at $t$ is $0$. Hence, for each $t$, the
probability that $X$ makes jump at $t$ is also $0$. Therefore we have
\[
Z^{\epsilon}(\tau)=0\text{ a.s. }%
\]

\item We remark that the equations $\left(  \ref{var1}\right)  -\left(
\ref{var2}\right)  $ are linear SDE and then by our assumptions on the
coefficients, they admit a unique solution.\newline
\end{enumerate}
\end{remark}

Let $\mathcal{R}$ denote the set of (Borel) measurable functions
$r:\mathbb{R}_{0}\rightarrow\mathbb{R}$ and define the Hamiltonian 
$H:\left[  0,T\right]  \times\mathbb{R}\times U\times\mathbb{R}%
\times\mathbb{R}\times\mathcal{R}\times\Omega\rightarrow\mathbb{R}$, to be%

\begin{align}
H(t,x,u,p,q,r)  &  :=H(t,x,u,p,q,\omega
)=f(t,x,u)+b(t,x,u)p\nonumber\\
&  \quad\quad+\sigma(t,x,u)q+%
{\textstyle\int_{\mathbb{R}_{0}}}
\gamma(t,x,u,\zeta)r(\zeta)\nu(d\zeta). \label{h}%
\end{align}
Let $(p^{\epsilon},q^{\epsilon},r^{\epsilon})\in\mathcal{S}^{2}\times
\mathbb{L}^{2}\times\mathbb{L}_{\nu}^{2}$ be the solution of the following
associated adjoint BSDE:
\begin{equation}%
\begin{cases}
dp^{\epsilon}(t) & =-\tfrac{\partial\tilde{H}}{\partial x}(t)dt+q^{\epsilon
}(t)dB(t)+%
{\textstyle\int_{\mathbb{R}_{0}}}
r^{\epsilon}(t,\zeta)\tilde{N}(dt,d\zeta);t\in\lbrack0,T],\\
p^{\epsilon}(T) & =\tfrac{\partial\tilde{g}}{\partial x}(\tilde{X}(T)),
\end{cases}
\label{p}%
\end{equation}
where
\[%
\begin{array}
[c]{ll}%
\tfrac{\partial\tilde{H}}{\partial x}(t) & =\tfrac{\partial\tilde{f}}{\partial
x}(t)+\tfrac{\partial\tilde{b}}{\partial x}(t)p^{\epsilon}(t)+\tfrac
{\partial\tilde{\sigma}}{\partial x}(t)q^{\epsilon}(t)+%
{\textstyle\int_{\mathbb{R}_{0}}}
\tfrac{\partial\tilde{\gamma}}{\partial x}(t,\zeta)r^{\epsilon}(t,\zeta
)\nu(d\zeta).
\end{array}
\]

\begin{lemma}
\label{est} The following holds,
\begin{equation}
Z^{\epsilon}(t)\rightarrow0\text{ as }\epsilon\rightarrow0^{+}\text{;}%
\quad\text{ for all }t\in\lbrack\tau,T]. \label{esz}%
\end{equation}%
\begin{equation}
(p^{\epsilon},q^{\epsilon},r^{\epsilon})\rightarrow(\hat{p},\hat{q},\hat
{r})\text{ when }\epsilon\rightarrow0^{+}, \label{esa}%
\end{equation}
where $(\hat{p},\hat{q},\hat{r})$ is the solution of the BSDE
\[%
\begin{cases}
d\hat{p}(t) & =-\tfrac{\partial\hat{H}}{\partial x}(t)dt+\hat{q}(t)dB(t)+%
{\textstyle\int_{\mathbb{R}_{0}}}
\hat{r}(t,\zeta)\tilde{N}(dt,d\zeta);t\in\lbrack0,T],\\
\hat{p}(T) & =\frac{\partial g}{\partial x}(\hat{X}(T)).
\end{cases}
\]

\end{lemma}

\noindent{Proof.} \quad By the Itô formula, we see that the solutions of the
equations $\left(  \ref{var1}\right)  -\left(  \ref{var2}\right)  ,$ are
\begin{align}
Z^{\epsilon}(t)  &  =Z^{\epsilon}(\tau+\epsilon)\exp(%
{\textstyle\int_{\tau+\epsilon}^{t}}
\{\tfrac{\partial\tilde{b}}{\partial x}(s)-\tfrac{1}{2}(\tfrac{\partial
\tilde{\sigma}}{\partial x}(s))^{2}+%
{\textstyle\int_{\mathbb{R}_{0}}}
[\log(1+\tfrac{\partial\tilde{\gamma}}{\partial x}(s,\zeta))-\tfrac
{\partial\tilde{\gamma}}{\partial x}(s,\zeta)]\nu(d\zeta)\}ds\nonumber\\
&  \quad\quad+%
{\textstyle\int_{\tau+\epsilon}^{t}}
\tfrac{\partial\tilde{\sigma}}{\partial x}(s)dB(s)+%
{\textstyle\int_{\tau+\epsilon}^{t}}
{\textstyle\int_{\mathbb{R}_{0}}}
\log(1+\tfrac{\partial\tilde{\gamma}}{\partial x}(s,\zeta))\tilde{N}%
(ds,d\zeta));\quad\tau+\epsilon\leq t\leq T. \label{z1}%
\end{align}
and%
\begin{equation}%
\begin{array}
[c]{ll}%
Z^{\epsilon}(t) & =\Upsilon(t)^{-1}[%
{\textstyle\int_{0}^{t}}
\Upsilon(s)(\tfrac{\partial\tilde{b}}{\partial u}(s)(u^{\epsilon}(s)-\hat
{u}(s))\\
& +{%
{\displaystyle\int_{\mathbb{R}_{0}}}
}\left(  \tfrac{1}{1+\tfrac{\partial\tilde{\gamma}}{\partial x}(s,\zeta
)}-1\right)  \tfrac{\partial\tilde{\gamma}}{\partial u}(s,\zeta)(v-\hat
{u}(s))\nu(d\zeta))ds+%
{\textstyle\int_{0}^{t}}
\Upsilon(s)\tfrac{\partial\tilde{\sigma}}{\partial u}(s)(v-\hat{u}(s))dB(s)\\
& +%
{\textstyle\int_{0}^{t}}
{\textstyle\int_{\mathbb{R}_{0}}}
\Upsilon(s)\left(  \tfrac{\tfrac{\partial\tilde{\gamma}}{\partial u}%
(s,\zeta)(v-\hat{u}(s))}{1+\tfrac{\partial\tilde{\gamma}}{\partial x}%
(s,\zeta)}-1\right)  \tilde{N}(ds,d\zeta)];\tau\leq t\leq\tau+\epsilon,
\end{array}
\label{z2}%
\end{equation}
where%
\[
\left\{
\begin{array}
[c]{ll}%
d\Upsilon(t) & =\Upsilon(t^{-})[-\tfrac{\partial\tilde{b}}{\partial
x}(t)+(\tfrac{\partial\tilde{\sigma}}{\partial x}(t)(u^{\epsilon}(t)-\hat
{u}(t)))^{2}\\
& +%
{\displaystyle\int_{\mathbb{R}_{0}}}
\left\{  \tfrac{1}{1+\tfrac{\partial\tilde{\gamma}}{\partial x}(t,\zeta
)}-1+\tfrac{\partial\tilde{\gamma}}{\partial x}(t,\zeta)\right\}  \nu
(d\zeta)dt-\tfrac{\partial\tilde{\sigma}}{\partial x}(t)dB(t)\\
& +%
{\displaystyle\int_{\mathbb{R}_{0}}}
\left(  \tfrac{1}{1+\tfrac{\partial\tilde{\gamma}}{\partial x}(t,\zeta
)}-1\right)  \tilde{N}(dt,d\zeta)];\tau\leq t\leq\tau+\epsilon,\\
\Upsilon(0) & =1.
\end{array}
\right.
\]
For more details see Appendix.\newline From \eqref{z2} we see that
$Z^{\epsilon}(\tau+\epsilon) \rightarrow0$ as $\epsilon\rightarrow0^{+}$, and
then from \eqref{z1} we deduce that $Z^{\epsilon}(t)\rightarrow0$ as
$\epsilon\rightarrow0^{+}$, for all $t$.\newline\linebreak The BSDE $\left(
\ref{p}\right)  $ is linear, and we can write the solution explicitly as
follows (see e.g. Theorem 2.7 in Øksendal and Sulem \cite{OS2}):
\begin{equation}%
\begin{array}
[c]{l}%
p^{\epsilon}(t)=\mathbb{E}[\tfrac{\Gamma(T)}{\Gamma(t)}\tfrac{\partial
\tilde{g}}{\partial x}(\tilde{X}(T))+%
{\textstyle\int_{t}^{T}}
\tfrac{\Gamma(s)}{\Gamma(t)}\tfrac{\partial\tilde{f}}{\partial x}%
(s)ds|\mathcal{F}_{t}];\qquad t\in\lbrack0,T],
\end{array}
\label{pfor}%
\end{equation}
where $\Gamma(t)\in\mathcal{S}^{2}$ is the solution of the linear SDE%
\[%
\begin{cases}
d\Gamma(t) & =\Gamma(t^{-})[\tfrac{\partial\tilde{b}}{\partial x}%
(t)dt+\tfrac{\partial\tilde{\sigma}}{\partial x}(t)dB(t)+%
{\textstyle\int_{\mathbb{R}_{0}}}
\tfrac{\partial\tilde{\gamma}}{\partial x}(t,\zeta)\tilde{N}(dt,d\zeta
)];\qquad t\in\lbrack0,T],\\
\Gamma(0) & =1.
\end{cases}
\]
From this, we deduce that $p^{\epsilon}(t)\rightarrow\hat{p}(t),$
$q^{\epsilon}(t)\rightarrow\hat{q}(t)$ and $r^{\epsilon}(t,\zeta
)\rightarrow\hat{r}(t,\zeta)$ as $\epsilon\rightarrow0^{+}.\newline\square
$\newline\linebreak We now state and prove the main result of this paper.

\begin{theorem}
[Necessary maximum principle]\label{Thm.Ness} Suppose $\hat{u}\in\mathcal{A}$
is maximizing the performance $\left(  \ref{perf}\right)  $. Then for all
$t\in\lbrack0,T)$ and all bounded $\mathcal{F}_{t}$-measurable $v\in V$, we
have%
\[
\tfrac{\partial H}{\partial u}(t,\hat{X}(t),\hat{u}(t))(v-\hat{u}(t))\leq0.
\]

\end{theorem}

\noindent{Proof.} \quad Consider
\begin{equation}
J(u^{\epsilon})-J(\hat{u})=I_{1}+I_{2}, \label{j}%
\end{equation}
where
\begin{equation}
I_{1}=\mathbb{E}[%
{\textstyle\int_{\tau}^{T}}
\{f(t,X^{\epsilon}(t),u^{\epsilon}(t))-f(t,\hat{X}(t),\hat{u}(t))\}dt],
\label{eq3.11}%
\end{equation}
and
\begin{equation}
I_{2}=\mathbb{E}[g(X^{\epsilon}(T))-g(\hat{X}(T))]. \label{eq3.12}%
\end{equation}
By the mean value theorem, we can write%

\begin{equation}
I_{1}=\mathbb{E[}%
{\textstyle\int_{\tau}^{\tau+\epsilon}}
\{\tfrac{\partial\tilde{f}}{\partial x}(t)Z^{\epsilon}(t)+\tfrac
{\partial\tilde{f}}{\partial u}(t)(u^{\epsilon}(t)-\hat{u}(t))\}dt+%
{\textstyle\int_{\tau+\epsilon}^{T}}
\tfrac{\partial\tilde{f}}{\partial x}(t)Z^{\epsilon}(t)dt], \label{I1}%
\end{equation}
and, applying the Itô formula to $p^{\epsilon}(t)Z^{\epsilon}(t)$ and by
$\left(  \ref{p}\right)  ,\left(  \ref{var1}\right)  $ and $\left(
\ref{var2}\right)  $, we have
\begin{align}
I_{2}  &  =\mathbb{E[}\tfrac{\partial\tilde{g}}{\partial x}(\tilde
{X}(T))Z^{\epsilon}(T)]=\mathbb{E[}p^{\epsilon}(T)Z^{\epsilon}(T)]\nonumber\\
&  =\mathbb{E[}p^{\epsilon}(\tau+\epsilon)Z^{\epsilon}(\tau+\epsilon
)]\nonumber\\
&  +\mathbb{E[}%
{\textstyle\int_{\tau+\epsilon}^{T}}
p^{\epsilon}(t)dZ^{\epsilon}(t)+%
{\textstyle\int_{\tau+\epsilon}^{T}}
Z^{\epsilon}(t)dp^{\epsilon}(t)+%
{\textstyle\int_{\tau+\epsilon}^{T}}
d\left\langle p^{\epsilon},Z^{\epsilon}\right\rangle (t)]\nonumber\\
&  =\mathbb{E[}p^{\epsilon}(\tau+\epsilon)(%
{\textstyle\int_{\tau}^{\tau+\epsilon}}
\{\tfrac{\partial\tilde{b}}{\partial x}(t)Z^{\epsilon}(t)+\tfrac
{\partial\tilde{b}}{\partial u}(t)(u^{\epsilon}(t)-\hat{u}(t))\}dt\nonumber\\
&  +%
{\textstyle\int_{\tau}^{\tau+\epsilon}}
\{\tfrac{\partial\tilde{\sigma}}{\partial x}(t)Z^{\epsilon}(t)+\tfrac
{\partial\tilde{\sigma}}{\partial u}(t)(u^{\epsilon}(t)-\hat{u}%
(t))\}dB(t)\nonumber\\
&  +%
{\textstyle\int_{\tau}^{\tau+\epsilon}}
{\textstyle\int_{\mathbb{R}_{0}}}
\{\tfrac{\partial\tilde{\gamma}}{\partial x}(t,\zeta)Z^{\epsilon}%
(t)+\tfrac{\partial\tilde{\gamma}}{\partial u}(t,\zeta)(u^{\epsilon}%
(t)-\hat{u}(t))\}\tilde{N}(dt,d\zeta))]\nonumber\\
&  +\mathbb{E[}%
{\textstyle\int_{\tau+\epsilon}^{T}}
\{p^{\epsilon}(t)\tfrac{\partial\tilde{b}}{\partial x}(t)Z^{\epsilon
}(t)-\tfrac{\partial\tilde{H}}{\partial x}(t)Z^{\epsilon}(t)+q^{\epsilon
}(t)\tfrac{\partial\tilde{\sigma}}{\partial x}(t)Z^{\epsilon}(t)\nonumber\\
&  +%
{\textstyle\int_{\mathbb{R}_{0}}}
r^{\epsilon}(t,\zeta)\tfrac{\partial\tilde{\gamma}}{\partial x}(t,\zeta
)Z^{\epsilon}(t)\nu(d\zeta)\}dt].
\end{align}
Using the generalized duality formula $\left(  \ref{geduB}\right)  $ and
$\left(  \ref{geduN}\right)  $, we get%
\begin{align}
I_{2}  &  =\mathbb{E[}%
{\textstyle\int_{\tau}^{\tau+\epsilon}}
\{p^{\epsilon}(\tau+\epsilon)(\tfrac{\partial\tilde{b}}{\partial
x}(t)Z^{\epsilon}(t)+\tfrac{\partial\tilde{b}}{\partial u}(t)(u^{\epsilon
}(t)-\hat{u}(t)))\nonumber\\
&  +\mathbb{E[}D_{t}p^{\epsilon}(\tau+\epsilon)|\mathcal{F}_{t}](\tfrac
{\partial\tilde{\sigma}}{\partial x}(t)Z^{\epsilon}(t)+\tfrac{\partial
\tilde{\sigma}}{\partial u}(t)(u^{\epsilon}(t)-\hat{u}(t)))\nonumber\\
&  +%
{\textstyle\int_{\mathbb{R}_{0}}}
\mathbb{E[}D_{t,\zeta}p^{\epsilon}(\tau+\epsilon)|\mathcal{F}_{t}%
]\{\tfrac{\partial\tilde{\gamma}}{\partial x}(t,\zeta)Z^{\epsilon}%
(t)+\tfrac{\partial\tilde{\gamma}}{\partial u}(t,\zeta)(u^{\epsilon}%
(t)-\hat{u}(t))\}\nu(d\zeta)\}dt]\nonumber\\
&  -\mathbb{E[}%
{\textstyle\int_{\tau+\epsilon}^{T}}
\tfrac{\partial\tilde{f}}{\partial x}(t)Z^{\epsilon}(t)dt], \label{I2}%
\end{align}
where by the definition of $H$ $\left(  \ref{h}\right)  $
\[%
\begin{array}
[c]{ll}%
\tfrac{\partial\tilde{f}}{\partial x}(t) & =\tfrac{\partial\tilde{H}}{\partial
x}(t)-\tfrac{\partial\tilde{b}}{\partial x}(t)p^{\epsilon}(t)-\tfrac
{\partial\tilde{\sigma}}{\partial x}(t)q^{\epsilon}(t)-%
{\textstyle\int_{\mathbb{R}_{0}}}
\tfrac{\partial\tilde{\gamma}}{\partial x}(t,\zeta)r^{\epsilon}(t,\zeta
)\nu(d\zeta).
\end{array}
\]
Summing $\left(  \ref{I1}\right)  $ and $\left(  \ref{I2}\right)  $, we obtain%
\begin{equation}%
\begin{array}
[c]{ll}%
I_{1}+I_{2} & =\mathbb{E[}%
{\textstyle\int_{\tau}^{\tau+\epsilon}}
\{\tfrac{\partial\tilde{f}}{\partial x}(t)+p^{\epsilon}(\tau+\epsilon
)\tfrac{\partial\tilde{b}}{\partial x}(t)+\mathbb{E[}D_{t}p^{\epsilon}%
(\tau+\epsilon)|\mathcal{F}_{t}]\tfrac{\partial\tilde{\sigma}}{\partial
x}(t)\\
& +%
{\textstyle\int_{\mathbb{R}_{0}}}
\mathbb{E[}D_{t,\zeta}p^{\epsilon}(\tau+\epsilon)|\mathcal{F}_{t}%
]\tfrac{\partial\tilde{\gamma}}{\partial x}(t,\zeta)\nu(d\zeta)\}Z^{\epsilon
}(t)dt]\\
& +\mathbb{E[}%
{\textstyle\int_{\tau}^{\tau+\epsilon}}
\{\tfrac{\partial\tilde{f}}{\partial u}(t)+p^{\epsilon}(\tau+\epsilon
)\tfrac{\partial\tilde{b}}{\partial u}(t)+\mathbb{E[}D_{t}p^{\epsilon}%
(\tau+\epsilon)|\mathcal{F}_{t}]\tfrac{\partial\tilde{\sigma}}{\partial
u}(t)\\
& +%
{\textstyle\int_{\mathbb{R}_{0}}}
\mathbb{E[}D_{t,\zeta}p^{\epsilon}(\tau+\epsilon)|\mathcal{F}_{t}%
]\tfrac{\partial\tilde{\gamma}}{\partial u}(t,\zeta)\nu(d\zeta)\}(u^{\epsilon
}(t)-\hat{u}(t))dt].
\end{array}
\label{sum}%
\end{equation}
By the estimate of $Z^{\epsilon}$ $\left(  \ref{esz}\right)  $, we get
\begin{equation}
\underset{\epsilon\rightarrow0^{+}}{\lim}X^{\epsilon}(t)=\hat{X}(t)\text{; for
all }t\in\lbrack\tau,T]\text{,} \label{estx}%
\end{equation}
and by $\left(  \ref{esa}\right)  $ we have
\begin{equation}
p^{\epsilon}(t)\rightarrow\hat{p}(t)\text{, }q^{\epsilon}(t)\rightarrow\hat
{q}(t)\text{ and }r^{\epsilon}(t,\zeta)\rightarrow\hat{r}(t,\zeta)\text{ when
}\epsilon\rightarrow0^{+}, \label{estpq}%
\end{equation}
where $(\hat{p},\hat{q},\hat{r})$ solves the BSDE
\begin{equation}
\left\{
\begin{array}
[c]{ll}%
d\hat{p}(t) & =-\tfrac{\partial\hat{H}}{\partial x}(t)dt+\hat{q}(t)dB(t)+%
{\textstyle\int_{\mathbb{R}_{0}}}
\hat{r}(t,\zeta)\tilde{N}(dt,d\zeta);\tau\leq t\leq T,\\
\hat{p}(T) & =\tfrac{\partial g}{\partial x}(\hat{X}(T)).
\end{array}
\right.  \label{adj}%
\end{equation}
Using the above and the assumption that $\hat{u}$ is optimal, we get%
\begin{align*}
0  &  \geq\underset{\epsilon\rightarrow0^{+}}{\lim}\tfrac{1}{\epsilon
}(J(u^{\epsilon})-J(\hat{u}))\\
&  =\mathbb{E[}\{\tfrac{\partial f}{\partial u}(\tau,\hat{X}(\tau),\hat
{u}(\tau))+\hat{p}(\tau)\tfrac{\partial b}{\partial u}(\tau,\hat{X}(\tau
),\hat{u}(\tau))+\mathbb{E[}D_{\tau}\hat{p}(\tau^{+})|\mathcal{F}_{t}%
]\tfrac{\partial\sigma}{\partial u}(\tau,\hat{X}(\tau),\hat{u}(\tau))\\
&  \text{ \ \ \ \ \ \ \ \ }+%
{\textstyle\int_{\mathbb{R}_{0}}}
\mathbb{E[}D_{\tau,\zeta}\hat{p}(\tau^{+})|\mathcal{F}_{t}]\tfrac
{\partial\gamma}{\partial u}(\tau,\hat{X}(\tau),\hat{u}(\tau),\zeta)\nu
(d\zeta)\}(v-\hat{u}(\tau))],
\end{align*}
where, by Theorem \ref{Th2.9},
\[%
\begin{array}
[c]{lll}%
\mathbb{E[}D_{\tau}\hat{p}(\tau^{+})|\mathcal{F}_{t}] & =\underset
{\epsilon\rightarrow0^{+}}{\lim}\mathbb{E[}D_{\tau}\hat{p}(\tau+\epsilon
)|\mathcal{F}_{t}] & =\hat{q}(\tau),\\
\mathbb{E[}D_{\tau,\zeta}\hat{p}(\tau^{+})|\mathcal{F}_{t}] & =\underset
{\epsilon\rightarrow0^{+}}{\lim}\mathbb{E[}D_{\tau,\zeta}\hat{p}(\tau
+\epsilon)|\mathcal{F}_{t}] & =\hat{r}(\tau,\zeta).
\end{array}
\]
Hence
\[
\mathbb{E}[\tfrac{\partial H}{\partial u}(\tau,\hat{X}(\tau),\hat{u}%
(\tau))(v-\hat{u}(\tau))]\leq0.
\]
Since this holds for all bounded $\mathcal{F}_{\tau}$-measurable $v$, we
conclude that%
\[
\tfrac{\partial H}{\partial u}(\tau,\hat{X}(\tau),\hat{u}(\tau))(v-\hat
{u}(\tau))\leq0\text{ for all }v.
\]
$\square$

\section{Linear-Quadratic Optimal Control with Constraints}

We now illustrate our main theorem by applying it to a linear-quadratic
stochastic control problem with a constraint, as follows:\newline Consider a
controlled SDE of the form
\[
\left\{
\begin{array}
[c]{ll}%
dX(t) & =u(t)dt+\sigma dB(t)+%
{\textstyle\int_{\mathbb{R}_{0}}}
\gamma(\zeta)\tilde{N}(dt,d\zeta);\quad t\in\lbrack0,T],\\
X(0) & =x_{0}\in\mathbb{R}.
\end{array}
\right.
\]
Here $u\in\mathcal{A}$ is our control process (see below) and $\sigma$ and
$\gamma$ is a given constant in $\mathbb{R}$ and function from $\mathbb{R}%
_{0}$ into $\mathbb{R}$, respectively, with
\[%
{\textstyle\int_{\mathbb{R}_{0}}}
\gamma^{2}(\zeta)\nu(d\zeta)<\infty.
\]
We want to control this system in such a way that we minimize its value at the
terminal time $T$ with a minimal average use of energy, measured by the
integral $\mathbb{E}[%
{\textstyle\int_{0}^{T}}
u^{2}(t)dt]$ and we are only allowed to use nonnegative controls. Thus we
consider the following constrained optimal control problem:

\begin{problem}
\label{pro} Find $\hat{u}\in\mathcal{A}$ (the set of admissible controls) such
that
\[
J(\hat{u})=sup_{u\in\mathcal{A}}J(u),
\]
where
\[
J(u)=\mathbb{E}[-\tfrac{1}{2}X^{2}(T)-\tfrac{1}{2}%
{\textstyle\int_{0}^{T}}
u^{2}(t)dt],
\]
and $\mathcal{A}$ is the set of predictable processes $u$ such that
$u(t)\geq0$ for all $t\in\lbrack0,T]$ and
\[
\mathbb{E}[%
{\textstyle\int_{0}^{T}}
u^{2}(t)dt]<\infty.
\]

\end{problem}

Thus in this case the set $V$ of admissible control values is given by
$V=[0,\infty)$ and we can use $U=V$. The Hamiltonian is given by
\[
H(t,x,u,p,q,r)=-\tfrac{1}{2}u^{2}+up+\sigma q+%
{\textstyle\int_{\mathbb{R}_{0}}}
\gamma(\zeta)r(\zeta)\nu(d\zeta),
\]
the adjoint BSDE for the optimal adjoint variables $\hat{p},\hat{q},\hat{r}$
is given by
\[
\left\{
\begin{array}
[c]{ll}%
d\hat{p}(t) & =\hat{q}(t)dB(t)+%
{\textstyle\int_{\mathbb{R}_{0}}}
\hat{r}(t,\zeta)\tilde{N}(dt,d\zeta);t\in\lbrack0,T],\\
\hat{p}(T) & =-\hat{X}(T).
\end{array}
\right.
\]
Hence
\begin{equation}
\hat{p}(t)= - \mathbb{E}[\widehat{X}(T) | \mathcal{F}_t].
\end{equation}

Theorem \ref{Thm.Ness} states that if $\hat{u}$ is optimal, then
\[
(-\hat{u}(t)+\hat{p}(t))(v-\hat{u}(t))\leq0;\quad\text{ for all }v\geq0.
\]
From this we deduce that
\[
\left\{
\begin{array}
[c]{cc}%
\text{ (i) if }\hat{u}(t)=0, & \text{ then }\hat{u}(t)\geq\hat{p}(t),\\
\text{ (ii) if }\hat{u}(t)>0, & \text{ then }\hat{u}(t)=\hat{p}(t).
\end{array}
\right.
\]
Thus we see that we always have $\hat{u}(t)\geq\max\{\hat{p}(t),0\}$. We claim
that in fact we have equality, i.e. that
\[
\hat{u}(t)=\max\{\hat{p}(t),0\} =\max \{- \mathbb{E}[\widehat{X}(T) | \mathcal{F}_t],0\}.
\]
To see this, suppose the opposite, namely that
\[
\hat{u}(t)>\max\{\hat{p}(t),0\}.
\]
Then in particular $\hat{u}(t)>0$, which by (ii) above implies that $\hat
{u}(t)=\hat{p}(t)$, a contradiction. We summarize what we have proved as follows:

\begin{theorem}
Suppose there is an optimal control $\hat{u}\in\mathcal{A}$ for Problem
\ref{pro}. Then
\[
\hat{u}(t)=\max\{\hat{p}(t),0\} =\max \{- \mathbb{E}[\widehat{X}(T) | \mathcal{F}_t],0\},
\]
where $(\hat{p},\hat{X})$ is the solution of the coupled forward-backward SDE
system given by
\begin{align*}
&
\begin{cases}
d\hat{X}(t) & =\max\{\hat{p}(t),0\}dt+\sigma dB(t)+%
{\textstyle\int_{\mathbb{R}_{0}}}
\gamma(\zeta)\tilde{N}(dt,d\zeta);\quad t\in\lbrack0,T],\\
\hat{X}(0) & =x_{0}\in\mathbb{R},
\end{cases}
\\
&
\begin{cases}
d\hat{p}(t) & =\hat{q}(t)dB(t)+%
{\textstyle\int_{\mathbb{R}_{0}}}
\hat{r}(t,\zeta)\tilde{N}(dt,d\zeta);t\in\lbrack0,T],\\
\hat{p}(T) & =-\hat{X}(T).
\end{cases}
\end{align*}

\end{theorem}

\begin{remark}
For comparison, in the case when there are no constraints on the control $u$,
we get from the well-known solution of the classical linear-quadratic control
problem (see e.g. Øksendal \cite{O}, Example 11.2.4) that the optimal control
$u^{\ast}$ is given in feedback form by
\[
u^{\ast}(t)=-\frac{X(t)}{T+1-t};\quad t\in\lbrack0,T].
\]

\end{remark}

\section{Appendix}

In this section, we give a solution of a general SDE with jumps. Let $X(t)$
satisfy the equation%
\[
\left\{
\begin{array}
[c]{ll}%
dX(t) & =(b_{0}(t)+b_{1}(t)X(t))dt+(\sigma_{0}(t)+\sigma_{1}(t)X(t))dB(t)\\
& +%
{\textstyle\int_{\mathbb{R}_{0}}}
(\gamma_{0}\left(  t,\zeta\right)  +\gamma_{1}\left(  t,\zeta\right)
X(t))\tilde{N}(dt,d\zeta)];t\in\left[  0,T\right]  ,\\
X(0) & =x_{0},
\end{array}
\right.
\]
for given $\mathbb{F}$-predictable processes $b_{0}(t),b_{1}(t),\sigma
_{0}(t),\sigma_{1}(t),\gamma_{0}\left(  t,\zeta\right)  ,\gamma_{1}\left(
t,\zeta\right)  $ with $\gamma_{i}\left(  t,\zeta\right)  \geq-1$ for
$i=0,1$.\newline Now suppose%
\begin{align*}
\Upsilon(t) &  =\exp[%
{\textstyle\int_{0}^{t}}
(-b_{1}(s)+\tfrac{1}{2}\sigma_{1}^{2}(s)-%
{\textstyle\int_{\mathbb{R}_{0}}}
\{\log(1+\gamma_{1}\left(  s,\zeta\right)  )-\gamma_{1}\left(  s,\zeta\right)
\}\nu(d\zeta))ds\\
&  \text{ \ \ \ \ \ \ \ \ \ \ \ \ \ \ \ \ }-%
{\textstyle\int_{0}^{t}}
\sigma_{1}(s)dB(s)+%
{\textstyle\int_{0}^{t}}
{\textstyle\int_{\mathbb{R}_{0}}}
\log(1+\gamma_{1}\left(  s,\zeta\right)  )\tilde{N}(ds,d\zeta)];t\in\left[
0,T\right]  .
\end{align*}
Then, $\Upsilon(t)=\exp(\Pi(t))$, where%
\[
\left\{
\begin{array}
[c]{ll}%
d\Pi(t) & =(-b_{1}(t)+\tfrac{1}{2}\sigma_{1}^{2}(t)-%
{\textstyle\int_{\mathbb{R}_{0}}}
\{\log(1+\gamma_{1}\left(  t,\zeta\right)  )-\gamma_{1}\left(  t,\zeta\right)
\}\nu(d\zeta))dt\\
& -\sigma_{1}(t)dB(t)+%
{\textstyle\int_{\mathbb{R}_{0}}}
\log(1+\gamma_{1}\left(  t,\zeta\right)  )\tilde{N}(dt,d\zeta);t\in\left[
0,T\right]  ,\\
\Pi(0) & =0.
\end{array}
\right.
\]
By the Itô formula as in Theorem 1.14 in Øksendal and Sulem \cite{OS1}, we
have%
\[
\left\{
\begin{array}
[c]{ll}%
d\Upsilon(t) & =\Upsilon(t^{-})[(-b_{1}(t)+\sigma_{1}^{2}(t)+%
{\textstyle\int_{\mathbb{R}_{0}}}
\{\tfrac{1}{1+\gamma_{1}\left(  t,\zeta\right)  }-1+\gamma_{1}\left(
t,\zeta\right)  \}\nu(d\zeta))dt\\
& -\sigma_{1}(t)dB(t)+%
{\textstyle\int_{\mathbb{R}_{0}}}
(\tfrac{1}{1+\gamma_{1}\left(  t,\zeta\right)  }-1)\tilde{N}(dt,d\zeta
)];t\in\left[  0,T\right]  ,\\
\Upsilon(0) & =1.
\end{array}
\right.
\]
Now put%
\[
Y(t)=X(t)\Upsilon(t).
\]
Then, again by the Itô formula, we obtain%
\begin{align}
dY(t) &  =d(X(t)\Upsilon(t))=X(t)d\Upsilon(t)+\Upsilon(t)dX(t)+d\langle
X,\Upsilon\rangle(t)\nonumber\\
&  =X(t)\Upsilon(t^{-})[(-b_{1}(t)+\sigma_{1}^{2}(t)-%
{\textstyle\int_{\mathbb{R}_{0}}}
\{\log(1+\gamma_{1}\left(  t,\zeta\right)  )-\gamma_{1}\left(  t,\zeta\right)
\}\nu(d\zeta)\nonumber\\
&  +%
{\textstyle\int_{\mathbb{R}_{0}}}
\{\tfrac{1}{1+\gamma_{1}\left(  t,\zeta\right)  }-1+\log(1+\gamma_{1}\left(
t,\zeta\right)  )\}\nu(d\zeta))dt\nonumber\\
&  -\sigma_{1}(t)dB(t)-%
{\textstyle\int_{\mathbb{R}_{0}}}
(\tfrac{1}{1+\gamma_{1}\left(  t,\zeta\right)  }-1)\tilde{N}(dt,d\zeta
)]\nonumber\\
&  +\Upsilon(t^{-})[(b_{0}(t)+b_{1}(t)X(t))dt+(\sigma_{0}(t)+\sigma
_{1}(t)X(t))dB(t)\nonumber\\
&  +%
{\textstyle\int_{\mathbb{R}_{0}}}
(\gamma_{0}\left(  t,\zeta\right)  +\gamma_{1}\left(  t,\zeta\right)
X(t))\tilde{N}(dt,d\zeta)]\nonumber\\
&  -(\sigma_{0}(t)+\sigma_{1}(t)X(t))\Upsilon(t^{-})\sigma_{1}(t)dt\nonumber\\
&  +{%
{\textstyle\int_{\mathbb{R}_{0}}}
}\Upsilon(t^{-})(\tfrac{1}{1+\gamma_{1}\left(  t,\zeta\right)  }-1)(\gamma
_{0}\left(  t,\zeta\right)  +\gamma_{1}\left(  t,\zeta\right)  X(t))\tilde
{N}(dt,d\zeta)\nonumber\\
&  +{%
{\textstyle\int_{\mathbb{R}_{0}}}
}\Upsilon(t^{-})(\tfrac{1}{1+\gamma_{1}\left(  t,\zeta\right)  }-1)(\gamma
_{0}\left(  t,\zeta\right)  +\gamma_{1}\left(  t,\zeta\right)  X(t))\nu
(d\zeta)dt.\label{dy}%
\end{align}
Rearranging terms, we end up with%
\begin{align*}
dY(t) &  =Y(t)[-%
{\textstyle\int_{\mathbb{R}_{0}}}
\{\log(1+\gamma_{1}\left(  t,\zeta\right)  )-\gamma_{1}\left(  t,\zeta\right)
+\tfrac{1}{1+\gamma_{1}\left(  t,\zeta\right)  }-1\\
&  +\log(1+\gamma_{1}\left(  t,\zeta\right)  )+(\tfrac{1}{1+\gamma_{1}\left(
t,\zeta\right)  }-1)\gamma_{1}\left(  t,\zeta\right)  \}\nu(d\zeta))dt\\
&  -%
{\textstyle\int_{\mathbb{R}_{0}}}
(\tfrac{1}{1+\gamma_{1}\left(  t,\zeta\right)  }-1)+\gamma_{1}\left(
t,\zeta\right)  +(\tfrac{1}{1+\gamma_{1}\left(  t,\zeta\right)  }-1)\gamma
_{1}\left(  t,\zeta\right)  \tilde{N}(dt,d\zeta)]\\
&  +\Upsilon(t^{-})[(b_{0}(t)-\sigma_{0}(t)\sigma_{1}(t)+{%
{\textstyle\int_{\mathbb{R}_{0}}}
}(\tfrac{1}{1+\gamma_{1}\left(  t,\zeta\right)  }-1)\gamma_{0}(t,\zeta
))\nu(d\zeta))dt\\
&  +\sigma_{0}(t)dB(t)+%
{\textstyle\int_{\mathbb{R}_{0}}}
(\gamma_{0}(t,\zeta)+(\tfrac{1}{1+\gamma_{1}\left(  t,\zeta\right)  }%
-1)\gamma_{0}\left(  t,\zeta\right)  )\tilde{N}(dt,d\zeta)].
\end{align*}
Consequently,%

\[
\left\{
\begin{array}
[c]{ll}%
dY(t)= & =\Upsilon(t^{-})[(b_{0}(t)+{%
{\textstyle\int_{\mathbb{R}_{0}}}
}(\tfrac{1}{1+\gamma_{1}\left(  t,\zeta\right)  }-1)\gamma_{0}(t,\zeta
))\nu(d\zeta))dt\\
& +\sigma_{0}(t)dB(t)+%
{\textstyle\int_{\mathbb{R}_{0}}}
\tfrac{\gamma_{0}(t,\zeta)}{1+\gamma_{1}\left(  t,\zeta\right)  }\tilde
{N}(dt,d\zeta)];t\in\left[  0,T\right]  ,\\
Y(0) & =x_{0}.
\end{array}
\right.
\]
Hence%
\begin{align*}
X(t)\Upsilon(t)  &  =Y(t)=y_{0}+%
{\textstyle\int_{0}^{t}}
\Upsilon(s)(b_{0}(s)+{%
{\textstyle\int_{\mathbb{R}_{0}}}
}(\tfrac{1}{1+\gamma_{1}\left(  s,\zeta\right)  }-1)\gamma_{0}(s,\zeta
))\nu(d\zeta))ds\\
&  +%
{\textstyle\int_{0}^{t}}
\Upsilon(s)\sigma_{0}(s)dB(s)+%
{\textstyle\int_{0}^{t}}
{\textstyle\int_{\mathbb{R}_{0}}}
\Upsilon(s)(\tfrac{\gamma_{0}(s,\zeta)}{1+\gamma_{1}\left(  s,\zeta\right)
})\tilde{N}(ds,d\zeta).
\end{align*}
Thus the unique solution $X(t)$ is given by
\begin{align*}
X(t)  &  =Y(t)\Upsilon(t)^{-1}=\Upsilon(t)^{-1}[x_{0}+%
{\textstyle\int_{0}^{t}}
\Upsilon(s)(b_{0}(s)+{%
{\textstyle\int_{\mathbb{R}_{0}}}
}(\tfrac{1}{1+\gamma_{1}\left(  s,\zeta\right)  }-1)\gamma_{0}(s,\zeta
)\nu(d\zeta))ds\\
&  +%
{\textstyle\int_{0}^{t}}
\Upsilon(s)\sigma_{0}(s)dB(s)+%
{\textstyle\int_{0}^{t}}
{\textstyle\int_{\mathbb{R}_{0}}}
\Upsilon(s)(\tfrac{\gamma_{0}(s,\zeta)}{1+\gamma_{1}\left(  s,\zeta\right)
})\tilde{N}(ds,d\zeta)];t\in\left[  0,T\right]  .
\end{align*}

\end{document}